\documentclass[12pt]{article}
\usepackage{latexsym,amssymb}
\def\bbbr{{\Bbb R}}     %Real numbers
\def\bbbp{{\Bbb P}}     %Projective space
\def\bbbc{{\Bbb C}}  % Complex numbers
\def\bbbz{{\Bbb  Z}}       %  numbers
\def\bbbn{{\Bbb  N}}       %  natural numbers
\def\bbbO{{\cal O}}     % sheaf O
\newtheorem{definition}{Definition}
\newtheorem{th1}{Theorem}[section]

\begin{document}
\title{Jacobians of singularized spectral curves
 and completely 
integrable systems}

\author{Lubomir Gavrilov\\
\normalsize \it Laboratoire Emile Picard, UMR 5580,
 Universit\'e Paul Sabatier\\
\normalsize \it 118, route de Narbonne,
 31062 Toulouse Cedex, France \\
gavrilov@picard.ups-tlse.fr\\ }
\date{April 2001}

\maketitle
\begin{abstract}
We state two recent results concerning the linearization of integrable 
systems on 
generalized Jacobians. Then we apply this to the (complexified) 
spherical pendulum.

\end{abstract}
\section{Introduction}

Let $M^J$ be the affine vector space of complex matrix polynomials $A(x)$ in a variable
$x$, of fixed degree $d$ and dimension $r$
$$
A(x)= J x^d+A_{d-1}x^{d-1}+...+A_0 , \; A_i \in {\bf gl}_r(\bbbc) 
$$
where 
$J\in {\bf gl}_r(\bbbc)$ is  a fixed 
matrix. The matricial polynomial Lax equations
\begin{equation}
\label{axa}
\frac{d}{dt} A(x)= [\frac{A^k(a)}{x-a},A(x)], k\in \bbbn, a\in \bbbc
\end{equation}
are well known to be Hamiltonian (with respect to several compatible Poisson 
structures on $M^J$)
and completely integrable. 
The corresponding Hamiltonian vector fields define
a complete set of commuting vector fields on the  isospectral manifolds
$$ M_P^J= \{A(x)\in M^J: det(A(x)-y I_r)= P(x,y) \} .$$ 
The corresponding 
complex flows are not complete, but may be completed by adding a suitable divisor
 $D_\infty$ called the Painlev\'e divisor of the system. The new manifold
$\bar{M}_P^J= M_P^J \cup D_\infty$ becomes a non-compact complex-analytic 
commutative group. The most remarkable fact concerning $\bar{M}_P^J$ is that 
it has
 a commutative algebraic group structure which is
 compatible with the structure of an analytic group. This algebraic group is just 
the generalized  Jacobian $J(X')$ of a suitable singularized spectral curve 
$X'$ (to be defined below). The variety
 $\bar{M}_P^J $ is non-compact and it admits several algebraic structures. 
The ``right" one is defined  by the symmetry 
group $G$ of the system (\ref{axa}). Namely, let
 $G= \bbbp{\bf GL}_r(\bbbc;J)$ be the subgroup of the projective group
$\bbbp{\bf GL}_r(\bbbc)$ formed by matrices which commute with $J$. In the 
applications $G$ is the symmetry group of the corresponding dynamical system (e.g. a rigid
body with an axis of symmetry has a symmetry group $\bbbc^*$ which is the 
complexified rotation group $S^1$).
The group $G$  acts on $M^J$ by conjugation,
 the
action is Poisson, and the reduced Hamiltonian system is completely integrable
 too. The action of  $G$ on $M_P^J$ can be extended on the completed variety
$\bar{M}_P^J$ and it is proper and free.
Therefore $\bar{M}_P^J$ can be considered as the
total space
 of
a holomorphic principal fiber bundle $\xi$
with base $\bar{M}_P^J /G$, structural group $G$, and natural projection map
$$
\bar{M}_P^J \stackrel{\phi}{\rightarrow}\bar{M}_P^J /G .
$$

The  fiber bundle $\xi$ is described as follows. 
When the spectral curve
$X$ defined by $\{(x,y)\in \bbbc^2: P(x,y)=0\}$ is smooth,
then the partially compactified variety $\bar{M}_P^J$ is smooth 
and bi-holomorphic to 
the generalized Jacobian variety
$J(X')$. 
The curve $X'$ is singular and as a topological space it is just $X$
with its ``infinite" points $\infty_1,\infty_2,...,\infty_r$ identified to a 
single point $\infty$. Thus $J(X')$ is a non-compact commutative algebraic
group and it can be
described as an extension of the usual Jacobian $J(X)$ by the algebraic group
$G=(\bbbc^*)^{s-1}\times \bbbc^{r-s}$, where $s\leq r$ is the number of distinct
eigenvalues of the leading term $J$
\begin{equation}
\label{aaa}
0\stackrel{ }{\rightarrow} G
 \stackrel{}{\rightarrow}
J(X')\stackrel{\phi}{\rightarrow}
J(X) \rightarrow 0 \; .
\end{equation}
As analytic spaces $J(X')$ and $J(X)$ are complex tori
$$
J(X')=\bbbc^{p_a}/\Lambda', \;J(X)=\bbbc^{p_g}/\Lambda
$$
where $\Lambda', \Lambda$ are lattices of rank $2p_g+s-1$ and $2p_g$ respectively,
$p_g$ is the genus of $X$, and $p_a=p_g+r-1$ is the arithmetic genus of $X'$.
The generalized Jacobian $J(X')$ can be also considered as 
the total space of
a holomorphic principal fiber bundle 
with base $J(X)$, projection  $\phi$, and structural group $G$. 
The group $G$ is then identified with the symmetry group 
$\bbbp{\bf GL}_r(\bbbc;J)$
of (\ref{axa}), and the  manifold $\bar{M}_P^J/G$ with  the usual 
Jacobian $J(X)=J(X')/G$. The algebraic 
description of the reduced invariant manifold $\bar{M}_P^J/G$ is a well
known  result proved by
 A.Beauville \cite{Beauville} and 
M.R.Adams, J.Harnad, J.Hurtubise \cite{Adams} (see also 
 M. Adler, P. van Moerbeke \cite{Adler}, and 
section 8.2. of the survey \cite{Reyman} by A.G.Reyman and 
M.A. Semenov-Tian-Shansky ).

The above can be illustrated on the following ``simple" example. 
Let $\tau _1,\tau _2 \in \bbbc$ be generic complex numbers. Consider the $\bbbz$-module
of rank three $\Lambda \subset \bbbc^2$
\begin{equation}
   \label{lattice}
	  \Lambda = \bbbz\left(
			  \begin{array}{c}
			     2 \pi i \\
			     0
			  \end{array}
			\right)
	  \oplus \bbbz \left(
			  \begin{array}{c}
			     0	   \\
			     2 \pi i
			  \end{array}
		     \right)
	  \oplus \bbbz
		     \left(
			 \begin{array}{c}
			    \tau_1  \\
			    \tau_2
			 \end{array}
		     \right) ,\ \
	  \;.
\end{equation}
 $\bbbc^2/\Lambda$ is a non--compact {\it algebraic group} and
it can be considered as a (non--trivial) extension of the
elliptic curve $\,\bbbc/\{ 2\pi i \bbbz \oplus \tau_1 \bbbz\}\,$
by $\, \bbbc^* \sim \bbbc/2 \pi i \bbbz$
\begin{equation}
   \label{suite}
		       0    \, \rightarrow \,
	 \bbbc   / 2\pi i \bbbz \, \rightarrow \,
	 \bbbc^2 / \Lambda    \, \stackrel{\phi}{\rightarrow} \,
	 \bbbc   / \{ 2\pi i \bbbz \oplus \tau_1 \bbbz\}
			    \, \rightarrow 0 \, ,
   \qquad
	 \phi(z_1,z_2)=z_1 \; .
\end{equation}
If  $d=2$, $r=2$ and $P(x,y)$ is a suitable polynomial, then it may be shown 
that 
$\bar{M}_P^J = \bbbc^2/\Lambda $, the symmetry group is 
$G=\bbbc/2 \pi i \bbbz$ and acts as
$$
G \times \bbbc^2/\Lambda \rightarrow \bbbc^2/\Lambda :
g\times 
\left(
			 \begin{array}{c}
			    z_1  \\
			    z_2
			 \end{array}
		     \right) 
\rightarrow \left(
			 \begin{array}{c}
			    z_1  \\
			    z_2 + g
			 \end{array}
		     \right) \; .
$$
The quotient group $\{\bbbc^2 / \Lambda \}/G$ is the elliptic curve
$\,\bbbc/\{ 2\pi i \bbbz \oplus \tau_1 \bbbz\}\,$ which is the Jacobian of the
spectral curve $X$ with an affine equation $\{ P(x,y)=0\}$. Finally
$\bbbc^2 / \Lambda $ is the generalized Jacobian of the singular spectral curve 
$X'$ which is obtained from $X$ by a suitable singularization 
( a procedure opposite to the
usual regularization of a singular curve). 
We note that $\bbbc^2 / \Lambda $ is also an extension of the elliptic curve
$\,\bbbc/\{ 2\pi i \bbbz \oplus \tau_2 \bbbz\}\,$ and hence it has a second 
algebraic structure. The ``right" algebraic structure is the one related to the symmetry
group $G=\bbbc/2 \pi i \bbbz$.

Generalized Jacobians as generic invariant manifolds of integrable systems 
appeared (implicitly) first in the papers of Jacobi 
\cite{Jacobi,Jacobi1,Jacobi2} which is 
easily seen from
the analytic expressions of the solutions (see also Klein and Sommerfeld 
\cite{Klein}).
This motivated the 
further 
study of the corresponding generalized Jacobi inversion problem (e.g. Clebsch an Gordan 
\cite{Clebsch}). The modern theory of generalized Jacobians, without relation to integrable 
systems, is due to Rosenlicht
\cite{Rosenlicht,Rosenlicht1}. In our 
exposition we shall follow Serre \cite{Serre}. In the modern literature on 
integrable systems generalized Jacobians appeared again in  
B.A. Dubrovin \cite{Dubrovin} and E. Previato \cite{Previato}. For more examples see
Fedorov \cite{Fedorov,Fedorov1,Alber} and Faye \cite{faye,faye1}.

In the present note we describe  the invariant manifold of the integrable system
(\ref{axa}) in the case when the leading term $J$ of the Lax matrix $A(x)$ is 
either
regular (Theorem \ref{gv}), or it is not regular but diagonalizable 
(Theorem \ref{viv}). 
The  proofs are given in 
\cite{Gav1,Vivolo}.
Then we apply Theorem \ref{gv} to the complexified spherical 
pendulum. The description which we obtain of the invariant manifold of the system
completes the recent results of M. Audin, F. Beukers and R. Cushman 
\cite{Audinp,Beukers,letter}. 

{\bf Acknowledgments}
I am grateful to 
 R. Cushman for sending me  the unpublished text \cite{letter} part of which is 
reproduced in the last section.

\section{ Singularized spectral curves and their Jacobians }
\label{general}

A polynomial
$$
P(x,y)= y^r+s_1(x)y^{r-1} + ... + s_r(x)
$$
 is called spectral, provided that 
the affine curve $\{(x,y)\in \bbbc^2: P(x,y)=0\}$ is the spectrum
of some polynomial $r \times r$ matrix $A(x)$
$$
P(x,y)= det(A(x) - y.I_r) \; .
$$
In this case $deg(s_i(x)) \leq i.d$, where $d$ is the degree of $A(x)$
\begin{equation}
\label{ax}
A(x)= A_d x^d+A_{d-1}x^{d-1}+...+A_0 , \; A_i \in {\bf gl}_r(\bbbc) \; .
\end{equation}
Consider the weighted projective space 
$\bbbp^2(d)=\bbbc^3 \backslash \{0\} / \bbbc^* $, where the $\bbbc^*$-action on
$\bbbc^3$ is defined by
$$
t\cdot (x,y,z) \rightarrow (tx,t^dy,tz), \; t\in \bbbc^* \; .
$$
$\bbbp^2(d)$ is a compact complex surface with one singular point $\{[0,1,0]\}=  \bbbp^2(d)_{sing}$.
 The affine curve  $\{(x,y)\in \bbbc^2: det(A(x) - y.I_r) = 0 \}$
is naturally embedded in $\bbbp^2(d)$, 
$$
\bbbc^2 \rightarrow \bbbp^2(d) : (x,y) \mapsto [x,y,1] ,
$$
and the condition $deg(s_i(x)) \leq i.d$ shows that its closure $X$ is contained 
in the smooth surface
$\bbbp^2(d)_{reg}=\bbbp^2(d)\backslash \{[0,1,0]\}$.
Let $x$ be an affine coordinate on $\bbbp^1$.
The surface $\bbbp^2(d)_{reg}$ is identified with the total space of the
holomorphic line bundle $\bbbO_{\bbbp^1}(d)$ with base $\bbbp^1$ and 
 projection 
$$
\pi :\bbbp^2(d)_{reg} \rightarrow \bbbp^1: [x,y,z] \rightarrow [x,z] .
$$
 The induced projection
\begin{equation}
\label{pi}
\pi : X \rightarrow \bbbp^1
\end{equation}
is a ramified covering of degree $r$, and over the affine plane $\bbbc$ it is simply the first projection
$$
\pi: \{(x,y)\in \bbbc^2 : P(x,y) = 0 \} \rightarrow \bbbc : (x,y) \rightarrow x \; .
$$
\begin{definition} (spectral curve of $A(x)$)
We define the spectral curve $X$ of the matrix polynomial $A(x)$
 (\ref{ax}),  to be the closure
of the affine curve $\{(x,y)\in \bbbc^2: det(A(x) - y.I_r) = 0 \}$ 
in the total space of the line bundle $\bbbO_{\bbbp^1}(d)$.
\end{definition}

From now on we fix the spectral polynomial $P(x,y)$ and  suppose that 
the spectral curve 
$X$ is {\it smooth} and irreducible.

We are going now to singularize the  curve $X$.
Let  $m = \sum_{i=1}^s n_i P_i$, 
$P_i \in X, n_i >0$, be an effective divisor on $X$. 
To the pair $(X, m)$  we associate a singular curve $X'= X_{reg}\cup \infty$, where
if  $S= \cup_{i=1}^s P_i$ is the support of $m$, then $X_{reg}= X-S$,  and $\infty$ is a single point. 
 The structure sheaf 
$\bbbO'$ of $X'\sim (X,m)$ is defined in the following way. Let $\bbbO_{X'}$ be the direct image of 
the structure sheaf $\bbbO=\bbbO_X$ under the canonical projection $X \rightarrow X'$.
Then
$$
\bbbO'_P = 
\left\{
\begin{array}{l}
\bbbO_P, \; \; P\in X_{reg} \\
\bbbc + i_\infty , \; \; P= \infty 
\end{array} 
\right.
$$
where  $i_\infty$ is the ideal of $\bbbO_\infty$ formed by the 
functions $f$ having a zero at $P_i$ of order at least $n_i$. 
Thus a regular function $f$ on $X'$ is a 
regular function $f$ on $X$, and such that for some $c\in \bbbc$ and any $i$ holds
$v_{P_i}(f-c) \geq n_i$, where $v_P(.)$ is the order function. If $p_g$ is the genus of $X$ then
the arithmetic genus $p_a$ of the singular curve $X'$ is $p_a=p_g+deg(m)-1$.

{\it Example} Let $m= P^++P^-$ be a divisor on the Riemann surface $X$. 
Then in a 
neighborhood of $\infty$ the singularized curve $X'$ is analytically isomorphic
either to the germ of analytical curve $xy=0$ ($P^+\neq P^-$), or to 
$y^2=x^3$ ($P^+=P^-$).

\begin{definition} (singularized spectral curve of $A(x)$)
\label{m}
If $\pi$ is the projection (\ref{pi}) and
$\infty =[1,0] \in \bbbp^1$ the ``infinite" divisor, then the effective divisor
$m=\pi^* (\infty) $
%= \sum_{i=1}^s n_i p_i$ 
is called 
modulus of the spectral curve $X$.
We have $deg(m)=r$ and we denote by $X'$  the singular curve associated to the regular curve $X$
and to the modulus $m$.
\end{definition}
The generalized Jacobian of the singular curve $X'$ is the
analytic manifold
$$
J(X')=  H^0(X,\Omega ^1( -m))^*/H_1(X_{reg},\bbbz) = \bbbc^{p_a}/ \Lambda' \; ,
$$
where  $\Lambda'$ is a rank $2p_g+s-1$ lattice, and $\Omega ^1(- m)$ is the sheaf
of meromorphic one-forms $\omega $, such that $(\omega ) \geq - m$. 
Similarly, for the usual Jacobian $J(X) \subset J(X') $, we have
$$
J(X)= H^0(X,\Omega ^1)^*/H_1(X,\bbbz) = \bbbc^{p_g}/ \Lambda \; ,
$$
where $\Lambda \subset \Lambda'$ is a rank $2p_g$ lattice. It may be shown that 
both $J(X')$ and $J(X)$ are commutative algebraic groups related as follows. 
$J(X')$ is a non-trivial {\it extension} of
$J(X)$ by the algebraic group $G = (\bbbc^*)^{s-1} \times \bbbc^{deg(m)-s}$
\begin{equation}
\label{extension}
0\stackrel{ }{\rightarrow} G
 \stackrel{ }{\rightarrow}
J(X')\stackrel{\phi}{\rightarrow}
J(X) \rightarrow 0
\end{equation}
where the map $\phi$ is induced by the natural homomorphisms
$$
H_1(X_{reg},\bbbz) \rightarrow H_1(X ,\bbbz),\;
 H^0(X,\Omega ^1) \rightarrow  H^0(X,\Omega ^1( -m)) \; .
$$
This means that the sequence (\ref{extension}) is exact in the usual sense
and moreover the algebraic structure of $G$ (respectively of $J(X)$) is induced
(respectively quotient) of the algebraic structure of $J(X')$.
 Note
that $J(X')$ is non-compact. Indeed, while the topological space of $J(X)$
is $(S^1)^{2p_g}$, the one of $J(X')$ is $(S^1)^{2p_g+s-1}\times \bbbr^{2deg(m)-s-1}$.

Let $M_P$ be the variety of $r\times r$ polynomial matrices of degree $d$ (\ref{ax}), which have a
fixed spectral polynomial $P(x,y)$
$$ M_P= \{A(x): det(A(x)-y I_r)= P(x,y) \} .$$
and let  $M_P^J=M_P \cap M^J$ be the isospectral manifold formed by matrices
 of the form (\ref{ax})
 with fixed leading term
$A_d=J$
\begin{equation}
\label{axj}
A(x)= J x^d+A_{d-1}x^{d-1}+...+A_0 , \; A_i \in {\bf gl}_r(\bbbc) \; .
\end{equation}
The stabilizer 
$$
\bbbp{\bf GL}_r(\bbbc;J) = \{ R \in \bbbp{\bf GL}_r(\bbbc): RJR^{-1}=J \} 
$$
 of 
$\bbbp{\bf GL}_r(\bbbc)$ at $J\in {\bf gl}_r(\bbbc)$ 
is a commutative algebraic group isomorphic to 
$(\bbbc^*)^{s-1} \times \bbbc^{deg(m)-s}$. It is a well known fact that
$M_P^J$ is a smooth manifold, $\bbbp{\bf GL}_r(\bbbc;J)$ 
acts freely and properly 
on $M_P^J$  by conjugation, and the quotient space 
$M_P^J/\bbbp{\bf GL}_r(\bbbc;J)$ 
is a smooth manifold  biholomorphic to $J(X)-\Theta$
\cite{Adams,Beauville}.

Consider the holomorphic principal fiber bundle $\xi$ with total space
$M_P^J$, structural group $\bbbp{\bf GL}_r(\bbbc;J)$,
base $M_P^J/\bbbp{\bf GL}_r(\bbbc;J)$, and natural projection map
$\varphi : M_P^J \rightarrow M_P^J/\bbbp{\bf GL}_r(\bbbc;J)$.
Consider also the associate principal bundle $\eta$ with base space $J(X)-\Theta$, total space
$J(X')-\Theta'$, projection map $\phi$, and structural group $G$ 
(see (\ref{extension})).

We are ready to formulate our first result
\begin{th1}(\cite[L. Gavrilov]{Gav1})
\label{gv}
Suppose that the spectral curve $X$ is smooth.
There exists a partial compactification of $M_P^J$ to a non-compact algebraic
manifold $\bar{M}_P^J= M_P^J \cup D_\infty$ which is isomorphic 
(as an algebraic manifold) to the generalized Jacobian $J(X')$ of the 
singularized spectral curve $X'$. The Hamiltonian flows (\ref{axa}) are 
translation invariant on $J(X')$. The action of the symmetry group 
$\bbbp{\bf GL}_r(\bbbc;J)$
on $M_P^J$ extends to an action on $\bar{M}_P^J$ which is free, proper and compatible
with the algebraic structure of $\bar{M}_P^J$. The principal algebraic bundles
$$0\rightarrow \bbbp{\bf GL}_r(\bbbc;J) \rightarrow
\bar{M}_P^J \rightarrow  \bar{M}_P^J / \bbbp{\bf GL}_r(\bbbc;J)
\rightarrow 0
$$
and
$$
0 \rightarrow G
 \rightarrow
J(X') \stackrel{\phi}{\rightarrow}
J(X) \rightarrow 0 \; 
$$
are isomorphic.
The Painlev\'e divisor $D_\infty= \bar{M}_P^J \backslash M_P^J$ is the pre-image 
$\Phi^{-1}(\Theta)$ of the Riemann theta divisor $\Theta \subset J(X)$. 
\end{th1}

A local analysis (at infinity) shows that the smoothness of the
compact spectral curve $X$  implies  the
regularity of the leading term $J$ of $A(x)$ (a matrix is  regular if 
its minimal and characteristic 
polynomials coincide). We shall suppose now that $J$ is not regular, and hence
 $X$ is never smooth (at infinity). We consider only the simplest case when
 $J$ is 
diagonal. Equivalently, if $m(\lambda )$ is the minimal polynomial of $J$, then
 $m(\lambda )=\Pi _{i=1}^k(\lambda - \lambda i)$, 
where $\lambda _i\neq \lambda _j$ for $i\neq j$, and
$k< r=dim J$. Geometrically this means that $X$ has only normal crossing 
singularities at infinity. Without loss of generality we suppose that $J$ is diagonal,
so $\{e_i\}_1^r$ are its eigenvectors. If $J$ is not regular the dimension of
the stabilizer $\bbbp{\bf GL}_r(\bbbc;J)$ is strictly bigger than $d-1$. 
This suggests that, in order to compensate the increase of the dimension of the stabilizer
of $J$, 
additional restrictions should be imposed on the subleading term $A_{d-1}$.
Let $\oplus_i E_i$ be a spectral decomposition of $\bbbc^r$, 
$E_i= Ker (J- \lambda _i I)$. To an endomorphism $u\in End(\bbbc^r)$ we associate
its restriction on $E_i$ denoted $u|_{E_i} \in End(E_i)$
$$
u|_{E_i} : v \rightarrow Pr_{E_i} u(v) \; .
$$
Let us fix now a basis $\{e_i\}_{i=1}^r$ of eigenvectors. 
To each matrix $A$ we associate by the above formula its restriction $A|_{E_i}$.
If $J$ is diagonal, then $A|_{E_i}$ is a ``diagonal" minor of $A$ of dimension
$dim E_i \times dim E_i$. To each proper subspace
 $E_i= Ker (J- \lambda _i I)$  of dimension $m_i>1$ we associate a diagonalizable regular
matrix $K_i$ of dimension $m_i$. Denote further the set of matrices $K_i$ by $K$.
 We replace
further the space $M^J$ by the affine space 
$$
M^{J,K}= \{A(x)= \sum_{i=1}^d A_i x^i: A_i \in {\bf gl}_r(\bbbc), A_d=J, 
A_{d-1}|_{E_{i}} = K_i \forall i \mbox{  such that  } dim E_i > 1  \; . \}
$$
and 
$M_P^J$ by
$$
M_P^{J,K} = \{ A(x) \in M^{J,K}, det(A(x) - yI) \equiv P(x,y) \}\; .
$$
The singularization $X'$ 
of the {\em singular} spectral curve $X$ 
is defined as follows.
Let $\tilde{X}$ be the desingularized $X$. Then apply Definition \ref{m} to the curve
$\tilde{X}$ and to the divisor $m$, where $m$ is the pre-image of 
the points at `infinity" on $X$.
It is seen that the topological space of $X$ is 
$X_{aff} \cup \infty$, where
$X_{aff}= \{(x,y)\in \bbbc^2: P(x,y)=0\}$. Finally we define the stabilizer 
(or the symmetry
group)
$$
\bbbp{\bf GL}_r(\bbbc;J,K)= \{R \in \bbbp{\bf GL}_r(\bbbc,J):
(RA_{d-1}R^{-1})|_{E_i}= K_i, \forall i \mbox{  such that  } dim E_i >1 \} \; .
$$

The analogue of Theorem \ref{gv} is the following
\begin{th1}(\cite[O.Vivolo]{Vivolo})
\label{viv}
Suppose that $J$ is a non-regular diagonal matrix , 
and let $K$ be a set of regular diagonal matrices
associated to each proper space of $J$ of dimension bigger than one. If the affine part
$X_{aff}$ of the spectral curve is smooth, then the variety $M_P^{J,K}$ 
is smooth.
There exists a partial compactification of $M_P^{J,K}$ to a non-compact algebraic
manifold $\bar{M}_P^{J,K}= M_P^{J,K} \cup D_\infty$ which is isomorphic 
(as an algebraic manifold) to the generalized Jacobian $J(X')$ of the 
singularized spectral curve $X'$. The Hamiltonian flows (\ref{axa}) are 
translation invariant on $J(X')$. The action of the symmetry group 
$\bbbp{\bf GL}_r(\bbbc;J,K)$
on $M_P^{J,K}$ extends to an action on $\bar{M}_P^J$ which is free, proper and compatible
with the algebraic structure of $\bar{M}_P^{J,K}$. The principal algebraic 
bundles
$$
0 \rightarrow \bbbp{\bf GL}_r(\bbbc;J,K) \rightarrow
\bar{M}_P^{J,K} \rightarrow  \bar{M}_P^{J,K} / \bbbp{\bf GL}_r(\bbbc;J,K)
\rightarrow 0
$$
and
$$
0\rightarrow G
\rightarrow
J(X')\stackrel{\phi}{\rightarrow}
J(X) \rightarrow 0 \; 
$$
are isomorphic.
The Painlev\'e divisor $D_\infty= bar{M}_P^{J,K} \backslash M_P^{J,K}$ is the pre-image 
$\Phi^{-1}(\Theta)$ of the Riemann theta divisor $\Theta \subset J(X)$. 
\end{th1}

See \cite{Vivolo} for applications.

\section{The  spherical pendulum}
The complexified spherical pendulum is described by the following system of
differential equations
\begin{equation}
\label{pendulum}
\begin{array}{rl}
\dot{x} & = \, v \\
 \dot{v} & = -e_3 + (\langle x, e_3 \rangle - \langle v,v \rangle )\, x
\end{array}
\end{equation}
constrained on the (invariant) manifold
\begin{displaymath}
M = \{ (x,v) \in {\bf C}^3 \times {\bf C}^3 \, \mid \, 
\langle x,x \rangle = 1, \, \, \langle x,v \rangle =0 \} \; .
\end{displaymath}
where $\langle x,y \rangle = \sum _{i=1}^3 x_i y_i$ and $e_1,e_2, e_3$ is the
standard basis of $\bbbc^3$. The algebraic manifold $M$ is the
tangent bundle $TQ$ of the complexified sphere $Q = \{x: \langle x,x\rangle = 1
\}$ and may be also identified to $T^*Q$ via the bilinear form $\langle
.,.\rangle$ on $T_xQ$. Therefore $TQ$ is a simplectic manifold and 
 with respect to this simplectic structure (\ref{pendulum})
is  a two degrees of freedom completely
integrable Hamiltonian system  with Hamiltonian function
$$
H(x,v)= \frac{1}{2} \langle v,v \rangle + x_3 \;.
$$
The second integral of
motion
$$
K(x,v)= x_1v_2-x_2v_1
$$
is associated  to the  symmetry group $G=SO(2,\bbbc)$ (rotations about the axes
$e_3$).
We denote the corresponding Hamiltonian vector fields by $X_H$, $X_K$.

Following Cushman \cite{letter}, consider the following automorphism of $M$
\begin{equation}
\label{automorphism}
\varphi: M \rightarrow M : (x,v) \rightarrow (y,u) = (x, x \times v) \; .
\end{equation}
It is straightforward to check that $\varphi^2(x,v)=(x,-v)$,
$\varphi^{-1}(y,u) 
= (y, u \times y)$.
The push forward under $\varphi $ of the Hamiltonian vector field
(\ref{pendulum})
is the Hamiltonian vector field $X_{\cal H}$ 
 corresponding to the Hamiltonian 
\begin{displaymath}
{\cal H}(y,u) = ({\varphi }^{-1})^{\ast }H(y,u) \, = \, 
\frac{1}{2} \langle u,u \rangle + \langle y, e_3 \rangle .
\end{displaymath}
The integral curves of $X_{\cal H}$ satisfy 
\begin{equation}
\begin{array}{rl}
\dot{y} & = \dot{x} \, = \, v \, = \, u \times y \\
 \dot{u} & = \dot{x} \times v + x \times \dot{v} \, 
= - x \times e_3 \, = \, 
e_3 \times y.
\end{array}
\label{eq-sec5two}
\end{equation}
Physically the above Hamiltonian system  is obtained 
by reducing the {\em right} $S^1$ symmetry given by rotating the 
Lagrange top about its symmetry axis at its momentum value $0$ 
and then choosing appropriate time and length scales. For more 
details see \cite[pp. 191--200]{Bates}.

As in \cite{Gav} we introduce the 
complex change of variables
\begin{equation}
\begin{tabular}{rlcrl}
$U_1$ & $=\, u_3 -i\, u_2 $ & \hspace{.25in} & $U_2$ & $=\, y_3 -i\, y_2$ \\
$V_1$ & $=\, u_1  $ & \hspace{.25in} & $V_2$ & $=\, y_1 $ \\
$W_1$ & $=\, u_3 +i\, u_2 $ & \hspace{.25in} & $W_2$ & $=\, y_3+i\, y_2 $. 
\end{tabular}
\label{change}
\end{equation}
Consider the matrices 
\begin{equation}
A(\lambda ) = \mbox{{\footnotesize $\left( \begin{array}{cc}
V_1 \lambda + V_2 & {\lambda }^2 + U_1\lambda +U_2 \\
 {\lambda }^2 + W_1 \lambda + W_2 & -(V_1 \lambda + V_2) 
\end{array} \right) $}} 
\label{eq-sec5five}
\end{equation}
and 
\begin{equation}
B(\lambda ) = \frac{A(\lambda ) - A(0)}{\lambda } \, = \, 
\mbox{{\footnotesize $\left( \begin{array}{cc} 
V_1 & \lambda + U_1 \\
 \lambda + W_1 & -V_1 \end{array} \right) $}}.
\label{eq-sec5six}
\end{equation}
\par
\noindent
\begin{th1}
(\cite[R. Cushman]{letter}). 

Equation (\ref{pendulum}) can be
written in Lax form as 
\begin{equation}
\label{laxpair}
2i\, \frac{d}{dt}A(\lambda ) = \left[ A(\lambda ), B(\lambda ) \right] \, = 
- \left[ A(\lambda ), \frac{A(0)}{\lambda } \right] . 
\label{eq-sec5seven}
\end{equation}
\end{th1}
The proof is straightforward. An equivalent Lax pair may be found in the
recent preprint of M. Audin (she attributes it to A. Reyman). A Lax
pair in the Lie algebra $so(3,1)$ appeared earlier in \cite{Reyman}.

A short calculation shows that characteristic polynomial of $A(\lambda )$ is
\begin{eqnarray}
\det ({\mu }\, I - A(\lambda )) &  = & 
{\mu }^2 + (V_1 \lambda + V_2)^2 + 
({\lambda }^2 + U_1\lambda + U_2)({\lambda }^2 + W_1\lambda + W_2) 
\nonumber \\
&  &\hspace{-1in} =\,  {\mu }^2 + {\lambda }^4 + (U_1+W_1){\lambda }^3 + 
(U_2 +W_2 +U_1W_1 +V^2_1){\lambda }^2 \nonumber \\
& & \hspace{-.5in} + \, (U_1W_2 +U_2W_1 +2\, V_1V_2) {\lambda } 
+ (U_2W_2 + V^2_2) .
\label{eq-sec5eight}
\end{eqnarray} 
Using  (\ref{change}) and (\ref{automorphism})  we find that 
\begin{eqnarray*}
U_1+W_1 & = & 2\, \langle u, e_3 \rangle \, = \,
 2\, \langle x\times v, e_3 \rangle  \, =\, 2k \\
U_2 +W_2 +U_1W_1 +V^2_1 & = & 2(\frac{1}{2} \langle u,u \rangle + 
\langle y,e_3 \rangle ) \\
& = &  2(\frac{1}{2} \langle v,v \rangle + 
\langle x,e_3 \rangle ) \, = \, 2h \\
U_1W_2 +U_2W_1 +2\, V_1V_2& = & 2\langle u,y \rangle \, = \, 2 \langle 
v \times x,x \rangle \, = \, 0 \\
U_2W_2 + V^2_2 & = & \langle y,y \rangle \, = \, \langle x,x \rangle \, = \, 
1.
\end{eqnarray*}
Hence we obtain the spectral curve 
\begin{equation}
X_{aff}= \{
0 = {\mu }^2 + F_c(\lambda ) \}
\label{eq-sec5nine}
\end{equation}
of the Lax equation (\ref{eq-sec5seven}) where 
\begin{equation}
F_c(\lambda ) = {\lambda }^4 +2k\, {\lambda }^3 + 2h\, {\lambda }^2 +1.
\label{eq-sec5ten}
\end{equation}
Let $X$ be the compactified and normalized curve $X_{aff}$. When $X_{aff}$ is
smooth, $X$ is an elliptic curve, $X=X_{aff}\cup\infty_1\cup\infty_2$. Let
$X'$ be the singularized curve $X$ with respect to the modulus
$\infty_1+\infty_2$ (Definition \ref{m}) (the topological space of $X'$ is
pictured in \cite[Fig.2]{Gav} and \cite[Fig.2]{Audinp}).  As a corollary of
Theorem
\ref{gv} we obtain the following
\begin{th1}
\label{thpendulum}
Suppose that the affine spectral   curve $X_{aff}$ of the spherical pendulum is
smooth. Then the invariant manifold $M_{hk}$ is also smooth.
There exists a partial compactification of $M_{hk}$ to a non-compact
algebraic manifold $\bar{M}_{hk}= M_{hk} \cup D_\infty$ which is isomorphic 
(as an algebraic manifold) to the generalized Jacobian $J(X')$ of the 
singularized spectral curve $X'$. The Hamiltonian flows $X_H$ and $X_K$ of the spherical
pendulum are 
translation invariant on $J(X')$. The action of the symmetry group 
$G=SO(2,\bbbc)=\bbbc^*$
on $M_{hk}$ extends to an action on $\bar{M}_{hk}$ which is free, proper and
compatible with the algebraic structure of $\bar{M}_{hk}$. The principal
algebraic bundles
$$
0\rightarrow \bbbc^* \rightarrow
\bar{M}_P^J \rightarrow  \bar{M}_{hk} / \bbbc^*
\rightarrow 0
$$
and
$$
0 \rightarrow \bbbc^*
 \rightarrow
J(X')\stackrel{\phi}{\rightarrow}
J(X) \rightarrow 0 \; 
$$
are isomorphic.
The Painlev\'e divisor $D_\infty$ is isomorphic to $\bbbc^*$
\end{th1}

{\bf Remarks}
The comlplexified spherical pendulum is a limiting case of the complexified
Lagrange top and the algebro-geometric structure of the latter was studied
in \cite{Gav}. 
The above theorem is the analogue of the main result of \cite{Gav}.
A big part of Theorem \ref{thpendulum} is proved by Beukers and Cushman
\cite{Beukers,letter}. It is shown in \cite{Beukers}  that the general
invariant manifold
$$
M_{hk}= \{(x,v)\in M: H(x,v)=h, K(x,v)=k \}
$$
is the total space of a principal $\bbbc^*$ bundle with base an elliptic
curve. The equations satisfied by the variables $U_i,V_i,W_i$ as well the Lax
pair  (\ref{laxpair}) are obtained by R. Cushman \cite{letter}.

Our second remark concerns the monodromy of the energy-momentum map
$$
M\rightarrow \bbbr^2: (x,v)\rightarrow (H(x,v),K(x,v))
$$
restricted to the set of its regular values. This is a ``classical" question
and we refer the reader to \cite{Bates}. As the regular fibers of the 
energy-momentum map are  real parts of the affine part of the generalized
Jacobian $J(X')$, then it suffices to study the monodromy of the flat bundle
$$
H_1(J(X'),\bbbz) \rightarrow (h,k)
$$
restricted to the set of real regular $(h,k)$ (and taking into account the real
structure of $J(X')$. But $H_1(J(X'),\bbbz)=H_1(X_{aff},\bbbz)$ and hence the
last question is equivalent to the study of the familiar Picard-Lefschetz
monodromy of the homology bundle associated to the fibration
$$
\{ ( \lambda,\mu ) \in \bbbc^2:
{\lambda }^4 +2k\, {\lambda }^3 + 2h\, {\lambda
}^2 +1 + \mu ^2=0 \} \rightarrow (h,k) \; .
$$
In the case of the Lagrange top this was (implicitly) used by O.Vivolo
\cite{Vivolo1,Vivolo2} and in the case of the spherical pendulum by M. Audin
\cite{Audinp}. 

The non-degeneracy of the frequency map is another remarkable feature
of the spherical pendulum. This has been first proved by E. Horozov (as conjecture
 of J.J. Duistermaat) in relation to applications of KAM theory \cite{horozov}.
Another proof, revealing the hidden algebraic geometry of the problem, 
can be found in \cite{gv}.


\begin{thebibliography}{99}
\bibitem{Adams}M.R.Adams, J.Harnad, J.Hurtubise, Isospectral Hamiltonian
flows in finite and infinite dimension, Commun.Math.Phys. 134, 555-585 (1990).
\bibitem{Adler} M.Adler, P. van Moerbeke, Linearization of Hamiltonian systems, Jacobi varieties and
representation theory, Advances in Math., vol.38, 318-379 (1980).
\bibitem{Alber} Alber, Mark S.; Fedorov, Yuri N. Wave solutions of evolution equations and Hamiltonian flows
on nonlinear subvarieties of generalized Jacobians. J. Phys. A 33 (2000), no. 47,
 8409--8425. 
\bibitem{Audinp} M. Audin, Monodromy of the spherical pendulum 
{\em via} algebraic geometry, preprint, Strasbourg (2001).
{\it Spinning Tops},  Cambridge Studies in Advanced Mathematics, vol.51,
Cambridge University Press, Cambridge, 1996.
\bibitem{Bates} R. Cushman, L. Bates, Global aspects of classical integrable systems,
Birkh\"auser, 1997.
\bibitem{Beauville} A. Beauville, Jacobiennes des courbes spectrales et syst\`emes hamiltoniens 
compl\`etement int\'egrables, Acta Math., vol.164, p. 211-235 (1990).
\bibitem{Beukers} F. Beukers, R. Cushman, The complex geometry of 
the spherical pendulum, preprint, Utrecht (2000).
\bibitem{Clebsch} A. Clebsch, P. Gordan, {\it Theorie der abelschen Functionen},
Leipzig, 1866.
\bibitem{letter} R. Cushman, Letter to L. Gavrilov dated 18 December 2000.
\bibitem{Dubrovin} B.A. Dubrovin, Geometry of Abelian varieties, Riemann surfaces and 
nonlinear equations, doctoral dissertation, Moscow, 1984 (in russian).
\bibitem{faye} I. Faye, Semi-abelian surfaces and integrable systems,
Proc. Edinburg Math. Soc.(2001) {\bf 44}, 249-265. 
\bibitem{faye1} I. Faye, Sur l'Int\'egrabilit\'e Alg\'ebrique et la Non-Int\'egrabilit\'e des Syst\`emes Hamiltoniens,
Ph.D. thesis, University of  Toulouse III, 2001.
\bibitem{Fedorov} Yu. Fedorov, 
Classical integrable systems and billiards related to generalized Jacobians. Acta
Appl. Math. 55 (1999), no. 3, 251--301.
\bibitem{Fedorov1}  Yu. Fedorov,
Discrete versions of some algebraic integrable systems related to generalized
Jacobians. SIDE III---symmetries and integrability of difference equations (Sabaudia, 1998), 147--160, CRM Proc.
Lecture Notes, 25, Amer. Math. Soc., Providence, RI, 2000. 
\bibitem{Gav} L. Gavrilov, A. Zhivkov, The complex geometry of Lagrange top, 
L'Enseignement Math\'ematique, t.44 (1998), p. 133-170.
\bibitem{Gav1} L. Gavrilov, Generalized Jacobians of spectral curves and
completely integrable systems, Math. Z. {\bf 230}, 487-508 (1999).
\bibitem{gv} L. Gavrilov, The real period function of the $A_3$ singularity and 
perturbations of the spherical pendulum,
Comp. Mathematica {\bf 123 } 167-184 (2000).
\bibitem{horozov}{E. Horozov, Perturbations of the spherical pendulum
and Abelian integrals, {\it J. Reine Angew. Math., 408, 114, (1990).}}
\bibitem{Jacobi} C. Jacobi, \"Uber eine neue Methode zur Integration der hyperelliptischen 
Differentialgleichungen
und \"uber die rationale Form ihrer vollst\"andigen algebraischen Integralgleichungen, Crelle, Bd.32,
p.220-226 (1846). 
\bibitem{Jacobi1} C. Jacobi, Sur la rotation d'un corps, Gesammelte Werke, Bd 2, 289-352, Chelsea, 1969.
\bibitem{Jacobi2} C. Jacobi, Fragments sur la rotation d'un corps tir\'es des manuscrits de
Jacobi et communiqu\'es par E. Lotner, Gesammelte Werke, Bd 2, 425-514, Chelsea, 1969.
\bibitem{Klein} F. Klein and A. Sommerfeld, Theorie des Kreisels, Teubner,
Leipzig, 1897-1910.
\bibitem{Moerbeke} P. van Moerbeke, D. Mumford, The spectrum of difference 
operators and algebraic curves, Acta Math., vol.143, p. 93-154 (1979).
\bibitem{Mumford} D.Mumford, {\em Tata Lectures on Theta II}, 
Progress in Mathematics, vol. 43, Birkh\"auser,
1984.
\bibitem{Previato} E.Previato, Hyperelliptic quasi--periodic and soliton 
solutions of the nonlinear Schr\"odinger
equation, Duke Math. J., vol.52, 329-377 (1985).
\bibitem{Reyman} A.G. Reyman, M.A. Semenov-Tian-Shansky,{\em Group 
theoretical
methods in the theory of finite dimensional Integrable systems}, in 
Dynamical Systems VII, Enciclopaedia
of Mathematical Sciences,{\bf 16} (Springer, 1994).
\bibitem{Rosenlicht} M. Rosenlicht,  Generalized Jacobian varieties, Ann. of Maths, {\bf 59}, 
1954, p.505-530.
\bibitem{Rosenlicht1}M. Rosenlicht, A universal mapping property of generalized 
Jacobian varieties,
Ann. of Maths, {\bf 66}, 
1957, p.80-88.
\bibitem{Serre} J.-P. Serre,{\it Groupes Alg\'ebriques et Corps de Classes},
Hermann, Paris, 1959.
\bibitem{Vivolo} O. Vivolo, Jacobians of singular spectral curves 
and completely integrable systems, Proc. Edinburg Math. Soc. 
{\bf 43}, 605-623 (2000).
\bibitem{Vivolo1}  O. Vivolo, The monodromy of  action variables of Lagrange top, 
preprint No 62, Laboratoire Emile Picard, Iniversit\'e de Toulouse III, 1995.
\bibitem{Vivolo2}  O. Vivolo, Syst\`emes Int\'egrables et Courbes Alg\'ebriques, 
Th\`ese, Universit\'e de Toulouse III, 1997.

\end{thebibliography}
\end{document}